\documentclass[12pt]{amsart}
\usepackage{graphicx}
\usepackage[latin1]{inputenc}
\usepackage{fullpage}
\usepackage{amsfonts,amsmath,amscd,amssymb}
\usepackage[all]{xy}
\usepackage{a4}

\newcommand{\Z}{{\mathbb{Z}}}
\newcommand{\Q}{\mathbb{Q}}
\newcommand{\cS}{{\mathcal{S}}}
\newcommand{\cM}{{\mathcal{M}^C}}

\newcommand{\pr}{\mathfrak{p}}
\newcommand{\qr}{\mathfrak{q}}

\newcommand{\mono}{\hookrightarrow}
\newcommand{\epim}{\twoheadrightarrow}

 \renewcommand{\Re}[2]{{\mathrm{Re}}(#1,#2)}
 \newcommand{\Hom}[2]{{\mathrm{Hom}}(#1,#2)}
 \newcommand{\Hm}{\mathrm{Hom}}
 \newcommand{\End}[1]{{\mathrm{End}}(#1)}
 \newcommand{\Ker}[1]{{\mathrm{Ke}\,}(#1)}
\newcommand{\Ke}{{\mathrm{Ke}\,}}
 
 \newcommand{\Jac}[1]{{\mathrm{Jac}}\,(#1)}
 \newcommand{\Soc}[1]{{\mathrm{Soc}}\,(#1)}
 
 \renewcommand{\Im}[1]{{\mathrm{Im}\,}(#1)}
\newcommand{\Ann}{{\mathrm{Ann}}}
\newcommand{\Qm}{{Q_{\mathrm{max}}}}

 \newcommand{\ZZ}{{\mathbb Z}}

\newcommand{\Tr}[2]{{\rm Tr}(#1,#2)}

\newcounter{zlist}

\theoremstyle{plain}

\newtheorem*{lem}{Lemma}
\newtheorem*{prop}{Proposition}
\newtheorem*{thm}{Theorem}

\newtheorem*{cor}{Corollary}

\theoremstyle{remark}
\newtheorem{defn}{Definition}

\title[Covering coalgebras and dual non-singularity]{Covering coalgebras and dual non-singularity}
\author{Christian Lomp}
\address{Departamento de Matem\'{a}tica Pura, Universidade do Porto, Portugal}
\author{Virg\'{\i}nia Rodrigues}
\address{Departamento de Matem\'{a}tica, Universidade Federal de Santa Catarina, Brasil}

\thanks{This work was carried out as part of the project {\it Interac\c{c}\~{o}es entre \'{a}lgebras e co-\'{a}lgebras} between the Universidade do Porto and Universidade Federal de Rio Grande do Sul and Universidade de S\~{a}o Paulo financed through GRICES (Portugal) and CAPES (Brasil) during the second author's visit to the UP. 
The first author was partially supported by Centro de Matem\'{a}tica da Universidade do Porto (CMUP), financed by FCT
(Portugal) through the programs POCTI (Programa Operacional Ci\^{e}encia, Tecnologia, Inova\c{c}\~{a}o) and POSI
(Programa Operacional Sociedade da Informa\c{c}\~{a}o), with national and European Community structural funds. He would also like to thank the {\em Funda\c{c}{\~a}o de Calouste Gulbenkian} for travel grant 78123, that allowed him to present this work at the {\em VIII Antalya Algebra Days}. The second author was supported by CAPES (Projeto 135/05 - BEX 1378/05-8) and would like to thank the Departamento de Matemática Pura for its hospitality.}

\begin{document}


\maketitle

\begin{abstract} Localisation is an important technique in ring theory and yields the construction of various rings of quotients. Colocalisation in comodule categories has been investigated by some authors where the colocalised coalgebra turned out to be a suitable subcoalgebra. Rather then aiming at a subcoalgebra we look at possible coalgebra covers $\pi:D\rightarrow C$ that could play the r\^ole of a dual quotient object. Codense covers will dualise dense (or rational) extensions; a maximal codense cover construction for coalgebras with projective covers is proposed.
We also look at a dual non-singularity concept for modules which turns out to be the comodule-theoretic property that  turns the dual algebra of a coalgebra into a non-singular ring. As a corollary we deduce that hereditary coalgebras and hence path coalgebras are non-singular in the above sense. We also look at coprime coalgebras and Hopf algebras which are non-singular as coalgebras.
\end{abstract}


\section{Introduction}

Embedding algebras into better ones where certain problems have
solutions is one of the major tools in ring theory. An analogous
tool for coalgebras does not always exists. Instead of embedding a
coalgebra into a better behaved  coalgebra one could also try to find a suitable better behaved coalgebra with a projection onto the first one - a covering coalgebra.

The maximal ring of quotients $\Qm(A)$ of an algebra $A$ is such an example of a universal object that has good properties in particular when the algebra in question is non-singular. Recall that an algebra $A$ is called left non-singular if left annihilators of non-zero elements are never essential as left ideal. This conditions is a kind of non-commutative torsion-freeness for $A$ and Johnson's theorem states that $A$ is left non-singular if and only if $\Qm(A)$ is von Neumann regular, i.e. the weak global dimension of $\Qm(A)$ is zero.

Throughout the text we will assume that rings $R$ are associative and have a unit. Furthermore we shall write homomorphisms of modules
opposite of scalars. A submodule $N$ of a left $R$-module $M$ is called {\it essential} ({\it small}) if  for all proper non-zero $L\subset M$: $N\cap L \neq 0$ ( $N+L\neq M$).
We denote a small submodule $N$ of $M$ by $N\ll M$. Given a module $M$ we denote by $\sigma[M]$ the category of submodules of factor modules of direct sums of copies of $M$ (see \cite{wisbauer}). For any pair of modules $X$ and $Y$ we denote the trace of $X$ in $Y$ by $\Tr{X}{Y} = \sum\{\mathrm{Im}(f) \mid f\in \Hom{X}{Y}\}$.

\subsection{The maximal ring of quotients} Given a ring $R$, an overring $S$ of $R$ is called a
{\it left ring of  quotients} if $\Hm_{R-}(S/R,S)=0$. The maximal
left ring of quotients $\Qm(R)$ of $R$ is any left ring of quotients
such that for any left ring of quotients $S$ of $R$ with embedding $j:R \mono S$  there exists
a unique ring homomorphism $\varphi:S\rightarrow \Qm(R)$ such that $j\varphi = \imath$ where $\imath:R\mono \Qm(R)$ denotes the embedding:

$$ \xymatrix{R \ar@{^{(}->}[rr]^j\ar@{_{(}->}[dr]^\imath & & S \ar[dl]^\varphi \\ & \Qm(R) }$$

The maximal left ring of quotients exists and can be constructed
as follows: Let $E=E(R)$ be the injective hull of $R$ as left
$R$-module. Then

$$ \Qm(R) := \{  x \in E \mid (x)f=0 \forall f \in \End{E} \textit{ with } (R)f=0 \}.$$

By construction $\Qm(R)$ is the submodule of $E$ that satisfies
$\Qm(R)/R=\Re{E/R}{E}.$ Where $\Re{X}{Y}=\bigcap \{ \Ker{f} \mid
f:X\rightarrow Y \}$ denotes the reject of $X$ in $Y$.

\subsection{Finite dimensional coalgebras}
Let us examine the following example:
Fix a field $k$ and consider the graph:
$$ \xymatrix{a \ar[r]^x & b}$$ and the path $k$-coalgebra $C$ associated to
this graph which has basis $\{a,x,b\}$ such that $a$ and $b$ are
group-like and $\epsilon(x)=0$ with
$$\Delta(x)=a\otimes x + x \otimes b.$$ Can we find a coalgebra
$D$ and a projection $\pi:D\rightarrow C$ such
that $D$ has better properties then $C$ ?
First note that $C^*$ is isomorphic to the path algebra of the
given graph and hence
$$ C^* = \left(  \begin{array}{cc}k & k \\ 0 & k\end{array}\right).$$
The maximal ring of quotients of $C^*$ is the full $2\times
2$-matrix ring $D=M_2(k)$ over $k$. Dualizing the inclusion $\imath:
C^* \rightarrow D$ we get a projection
$$D^* \rightarrow C^{**} \simeq C$$
where $D^*$ is the $2 \times 2$-matrix coalgebra.

\subsection{} Analogous to the example above, we can always choose $D$ to be the dual coalgebra of the maximal ring of quotient of a finite dimensional coalgebra in order to obtain a suitable coalgebra cover, as the following theorem shows:
\begin{thm}\label{fdCoalgebras_maxQuoring} Let  $C$ be  a finite dimensional $k$-coalgebra,
then $D=\left(Q^{r}_{max}(C^*)\right)^*$  is a finite dimensional coalgebra and there exists
a surjective coalgebra homomorphims
$\pi: D \rightarrow C$ whose kernel is small as a right $C$-subcomodule of $D$.
\end{thm}

\begin{proof}
Since $C$ is finite dimensional, it is a left and right semiperfect coalgebra.
Let $P$ be a projective cover of $C$ as right $C$-comodule with epimorphism $\pi:P\rightarrow C$.
Since $C$ is finitely generated as left $C^*$-module, $P$ is also finitely generated as left $C^*$-module and
hence finite dimensional. Since $P$ is a projective right $C$-comodule, $P^*$ is an injective right $C^*$-module
(by \cite[9.5]{BrzezinskiWisbauer}).
Moreover as $\pi^*:C^* \rightarrow P^*$ is an essential embedding,
$P^*$ is isomorphic to the injective hull $E(C^*)$ of $C^*$ as right $C^*$-module.
Since $Q^{r}_{max}(C^*) \subseteq E(C^*)$, it is also finite dimensional.
Hence $D=\left(Q^{r}_{max}(C^*)\right)^*$ is a finite dimensional coalgebra and the transpose $\imath^*:D\epim C$ of the
algebra embedding $\imath:C^* \rightarrow Q^{r}_{max}(C^*)$ is a surjective coalgebra homomorphism. Since $\imath$ is an essential monomorphism, $\pi$ is a small epimorphism.
\end{proof}

\subsection{} Let $K$ be a field and $\Gamma$ be a quiver, i.e. a directed graph with finitely many vertices $\Gamma_0$ and finitely many arrows $\Gamma_1$ and without cycles.
The path $K$-coalgebra $C$ associated to $\Gamma$ is the vector space whose basis are all paths in $\Gamma$ and with comultiplication $\Delta(w)=\sum_{uv=w} u\otimes v$. For each vertex $i\in \Gamma_0$ denote  by $v_i$ the unique path of length zero at vertex $i$. Note that $C$ is finite dimensional and $C^*$ is isomorphic to the path algebra associated to $\Gamma$. Since by \cite[13.25]{Lam} the right maximal ring of quotients of a right artinian right non-singular ring $A$ is isomorphic $\End{\Soc{A_A}}$, we only need to determine the right socle of $C^*$ to describe the right maximal ring of quotients of $C^*$.
Let $A$ be the path algebra associated to $\Gamma$. Denote by $\Gamma_{sink}$ the set of {\it terminal} vertices $i\in \Gamma_0$, i.e. those vertices from where no arrow starts. Note that for any $i\in \Gamma_{sink}: v_iA=v_iK$ is a minimal right ideal of $A$. Moreover for any path $p$ in $A$ which ends at a terminal vertex $i$, the cyclic right ideal $pA$ is  a minimal right ideal and isomorpic to $v_iA$ since both have the same maximal right ideal $M_i$ generated by all paths except $v_i$. On the other hand let $I$ be a minimal right ideal of $A$, then $I=\gamma A$ for some linear combination $\gamma=\sum_{j=1}^n \lambda_j p_j$ of distinct paths $p_j$ and  non-zero coeficients $\lambda_j$. Let $i'$ be the vertex where the path $p_1$ ends and choose a path $q$ from $i'$ to some terminal vertex $i$. Then 
$I=\gamma q A$, since $I$ was minimal. Note that $qM_i=0$ implies that $IM_i=0$, i.e. the annihilator of $I$ is the maximal right ideal $M_i$. Hence $I\simeq v_iA$. Moreover $\gamma q$  can be written as a linear combination of paths ending at $i$, i.e. $\gamma q = \sum \lambda_j p'_j$ where all paths $p'_j$ end at $i$. Hence $I\subseteq \bigoplus p'_jA$.
For any terminal vertex $i\in \Gamma_{sink}$ denote by $P_i$ the set of paths ending at $i$ and set $n_i=|P_i|$. Then we just showed that
$$ \Soc{A_A} = \bigoplus_{i\in \Gamma_{sink}} \left( \bigoplus_{p\in P_i} pA\right) \simeq \bigoplus_{i\in\Gamma_{sink}} (v_iA)^{n_i}.$$
By \cite{Lam} the maximal right ring of quotients of $A$ is isomorphic to the endomorphism ring of $\Soc{A_A}$:
$$Q_{max}^r(A) \simeq \End{\Soc{A_A}} \simeq \prod_{i\in \Gamma_{sink}} \End{(v_iA)^{n_i}} \simeq 
\prod_{i\in \Gamma_{sink}} M_{n_i}(K),$$
where $M_n(K)$ denotes the ring of $n\times n$-matrizes over $K$.

Going back to our path coalgebra we have now a projection of coalgebras of a direct product of matrix coalgebra onto $C$, i.e.
$$ \prod_{i\in \Gamma_{sink}} M_{n_i}^c(K) \epim C.$$
Here $M_n^c(K)={\left(M_n(K)\right)}^*$ denotes the $n\times n$-matrix coalgebra with basis $\{E_{ij}\}_{1\leq i,j\leq n}$, comultiplication
$$ \Delta(E_{ij})=\sum_{l=1}^n E_{il}\otimes E_{lj}$$ and counit
$\epsilon(E_{ij})=\delta_{i,j}$.  

\subsection{}\label{dividedPower}
In case of an infinite dimensional path coalgebra, how can we
obtain a covering coalgebra like the matrix coalgebra in our
example ? For instance for the divided power coalgebra, that is the path coalgebra associated to the graph

$$\xygraph{. :@(ul,ur)^x "."}$$

We will see that there is no apropriate coalgebra cover in the sense defined below.

%



\section{A module-theoretic approach to covering coalgebras}
A module extension $X\mono Y$ is called {\it dense} if
$\Hm(Z/X,Y)=0$ for all $X\subset Z \subset Y$.
In \cite{FindlayLambek} Findlay and Lambek proved that the maximal ring of quotient
$Q$ of a ring $R$ is the maximal dense extension of $R$ in the category of $R$-modules.
We will give a module theoretic approach in covering coalgebras using {\it codense covers} of modules:

\subsection{Codense covers of modules}

A module $Y$ is called a {\it cover} of $X$ if there exists an
epimorphism $\pi:Y\epim X$. The cover $Y$ is said to be {\it small} if
$\Ker{\pi} \ll Y$ and a cover $Y$ is called a {\it codense cover} of
$X$ if  $\Ker{\pi}$ is a {\it codense submodule} of $Y$, that is
$\Hom{Y}{\Ker{\pi}/L}=0$ for all $L\subseteq \Ker{\pi}$. As a dualisation of
dense extensions, codense covers were introduced by Courter in
\cite{Courter} where they are called {\it co-rational extensions}.
Since the term {\it rational module} has a different meaning in the
coalgebraic setting, we prefer to refer to 'dense extensions' and
'codense covers'  instead. A non-trivial example of
a codense cover is the projection $\Q \epim \Q/\Z$, which is codense
since $\Hom{\Q}{\Z/n\Z}=0$  for all $n$.

\subsection{} Some properties of codense covers can be easily checked:

\begin{lem}\label{prelema} \label{projectiveCoverCodense}  Let $Z$ be a cover of $X$ in $\sigma[M]$.
\begin{enumerate}
\item If $Z$ is a codense cover of $X$, then it is a small cover.
\item If $Z$ is a small cover of $X$ and $\pi:X\epim Y$ is a codense cover then $\Hom{Z}{\Ker{\pi}}=0$.
\item If $Z$ is a projective cover of $X$ in $\sigma[M]$ then a cover $\pi:X \epim Y$ is codense if and only if $\Hom{Z}{\Ker{\pi}}=0$.
\end{enumerate}
\end{lem}


\begin{proof} (1) Let $\pi:Z\epim X$ be a codense cover. Suppose $\Ker{\pi}+Y=Z$, then the canonical projection
$$ Z \rightarrow Z/Y \simeq \Ker{\pi}/(\Ker{\pi}\cap Y)$$ is zero by hypothesis. Thus $Z=Y$ and $\Ker{\pi}\ll Z$.\\

(2) Let $p:Z\rightarrow X$ be a small epimorphism and $f\in \Hom{Z}{\Ker{\pi}}$.
Extending $f$ to an homomorphism $$g:X=Z/\Ker{p}\rightarrow \Ker{\pi}/(\Ker{p})f,$$
mapping $z+\Ker{p} \mapsto (z)f + (\Ker{\pi})f$, we have $g=0$ since $X$ is a codense cover of $Y$. Thus $(\Ker{p})f=\Im{f}$. But as $(\Ker{p})f \ll \Im{f}$, we must have $\Im{f}=0$, i.e. $f=0$.\\

(3) Since $Z$ is projective cover of $X$ there exists a small epimorphism
 $p:Z\epim X$. For any $U\subseteq \Ker{\pi}$ and $f:X\rightarrow \Ker{\pi}/U$ we have $pf:Z\rightarrow \Ker{\pi}/U$. Since $Z$ is projective there exists $g:Z\rightarrow \Ker{\pi}$ which is zero by hypothesis. Hence $pf=0$ and $f=0$ as $p$ is an epimorphism.
\end{proof}

\subsection{}
Dual to the definition of a  maximal dense extension of a module, we define a maximal codense cover as follows:

\begin{defn} Let $X,Y \in \sigma[M]$. A codense cover $p:Y\epim X$ is
called a {\it maximal codense cover in $\sigma[M]$ } if for any codense cover
$\pi:Z\epim X$ there exists a unique epimorphism $\psi:Y\rightarrow Z$ such that $\psi \pi = p$.
$$ \xymatrix{Z \ar@{->>}[rr]^{\pi} & & X  \\ & Y\ar[ul]^\psi \ar@{->>}[ur]^p }$$
\end{defn}

Note that our definition differs from Courter's in \cite{Courter}.

\subsection{} As it was to expect, in case projective covers exists a dual construction like Findlay and Lambek's allows to construct a maximal codense cover for modules:

\begin{thm}\label{MaxmialCodenseCovers} Let $X\in \sigma[M]$ have a  projective cover $P$ in $\sigma[M]$.
Denote by $\pi:P\rightarrow X$ the projection and $T:=\Tr{P}{\Ker{\pi}}$.
Then $\widetilde{P}=P/T$ is a maximal codense cover of $X$ in $\sigma[M]$ with induced epimorphism
$\widetilde{\pi}:\widetilde{P}\rightarrow X$.
\end{thm}

\begin{proof} Note that $\Ker{\widetilde{\pi}} = \Ker{\pi}/T$ and as $P$ is a projective cover of $\widetilde{P}$, $\Hm(P,\Ker{\pi}/T)=0$. By Lemma \ref{projectiveCoverCodense} $\widetilde{\pi}:\widetilde{P}\epim X$ is a codense cover.
Let $p:Z\epim X$ be any other codense cover of $X$ in $\sigma[M]$.
By the projectivity of $P$ there exist $\psi:P\rightarrow Z$ such that $\psi p = \pi$.
As $(T)\psi p = (T)\pi = 0$ we deduce $$ P\Hom{P}{\Ker{\pi}}\psi = T\psi \subseteq \Ker{p}.$$
Since by Lemma \ref{projectiveCoverCodense} $\Hom{P}{\Ker{p}}=0$, $(T)\psi = 0$.
Hence $\psi$ lifts to a homomorphism $\widetilde{\psi}:\widetilde{P}\rightarrow Z$ with $\widetilde{\psi} p = \widetilde{\pi}$.

$\widetilde{\psi}$ is unique because if there existed another map $\phi:\widetilde{P}\rightarrow Z$ with $\phi p = \widetilde{\pi}$, then $\psi-\phi\in \Hom{P}{\Ker{p}}=0$ by Lemma \ref{projectiveCoverCodense} (here we consider $\phi$ as a map from $P$ to $Z$).
\end{proof}





\subsection{}
For a finite dimensional coalgebra $C$ we saw in \ref{fdCoalgebras_maxQuoring} that $D=\left(Q^{r}_{max}(C^*)\right)^*$ is a small cover of $C$. Actually as it was to expect, $D$ is a maximal codense cover of $C$ in the category of right $C$-comodules:

\begin{thm}
Let $C$ be a finite dimensional coalgebra over a field $k$, then $D=\left(Q^{r}_{max}(C^*)\right)^*$ is a maximal codense cover of $C$ in $\cM$.
\end{thm}

\begin{proof}
By transposing the embedding $\imath: C^*\mono Q^{r}_{max}(C^*)$ we obtained a small cover
$\pi:D\epim C$ in \ref{fdCoalgebras_maxQuoring}, where $D=\left(Q^{r}_{max}(C^*)\right)^*$ and $\pi = \imath^*$.
The kernel $K$ of $\pi$ is isomorphic to $\left(Q^{r}_{max}(C^*)/C^*\right)^*$.
Note that the dual of any factor comodule $K\epim L$ is a right $C^*$-submodule
$L^*$ of $Q^{r}_{max}(C^*)/C^*$. Hence the transpose map of any right $C$-colinear map
$g:C\rightarrow K/L$ yields a right $C^*$-linear  map $g^*:(K/L)^* \rightarrow C^*$ which could be extended
to a right $C^*$-linear map from $Q^{r}_{max}(C^*)/C^*$ to $E(C^*)$ and must be zero (where $E(C^*)$ denotes the injective
hull of $C^*$ as right $C^*$-module).
Hence $D$ is a codense cover of $C$. The maximality follows now by a similar argument, taking into account
that any codense cover $D'$ of $C$ in $\cM$ would be finitely generated as comodule and hence finite dimensional.
\end{proof}

\subsection{}
We will now turn to some examples of modules that are equal its own
maximal codense cover.
The next Lemma is probably known, but we were
unable to find a reference:

\begin{lem}\label{uniserialInjectivesPID}
Every indecomposable non-faithful injective module over a principal ideal domain is uniserial.
\end{lem}

\begin{proof} Let $D$ be a principal ideal domain and $M$ an indecomposable non-faithful injective $D$-module.
By Matlis Theorem \cite{Matlis} $M=E(D/\pr)$ for some non-zero prime ideal $\pr=Dp$ of $D$.
Since $D$ is a Dedekind domain, the localisation of $D$ by $\pr$: $D_\pr$ is a discrete valuation ring.
Hence $D_\pr, Q$ and $Q/D_\pr$ are uniserial $D_\pr$-modules. Take any $D$-submodule $N \subseteq Q/D_\pr$.
We will show that $N$ is also a $D_\pr$-module. For any $a\not\in \pr=Dp$ and $n=x/y + D_\pr \in N$ with
 $y=up^k \in \pr$ and $p\nmid u$.
Hence $1=ra+sp^k$ for some $r,s \in D$. This implies that $\frac{1}{a} -r = \frac{sp^k}{a} \in D_\pr$. Therefore
$$\frac{1}{a}n - rn = \frac{sp^k}{a}\frac{x}{up^k}  = \frac{sx}{au} \in  D_\pr \Rightarrow \frac{1}{a}n = rn + D_\pr.$$
Hence the action of $1/a$ on an element $n$ in $Q/D_\pr$ is given by a $D$-scalar multiplication.
This shows that $Q/D_\pr$ is a uniserial $D$-module. Since $Q/D_\pr$ is injective and contains a simple
$D_\pr$-submodule which is isomorphic to $D_\pr/\pr D_\pr \simeq D/\pr$, we have
that $$M\simeq E(D/\pr) \simeq Q/D_\pr$$
is a uniserial $D$-module. Note that all its submodules are of the form $D/\pr^i$.
\end{proof}

\subsection{}
The next theorem states that indecomposable injectives over suitable rings do not have proper codense covers and as we will see below applies in particular to the case of the divided power coalgebra mentioned in \ref{dividedPower}.
A module $M$ is called {\it couniform} or {\it hollow} if every proper submodule is small.

\begin{thm}\label{PID_small_covers}
The only possible small covers of a  non-faithful indecomposable
injective module $M$ over a principal ideal domain $D$ are $M$ and
the quotient field $Q$ of $D$.
\end{thm}

\begin{proof} By a theorem of Matlis \cite[Prop 3.1]{Matlis} $M=E(D/\pr)$ for some maximal ideal $\pr$.
Furthermore $M$ is uniserial by \ref{uniserialInjectivesPID}. Let $\pi: P\rightarrow M$ be
a small cover. Then $P$ is hollow, since $M$ is uniserial and whenever $P=D + E$,
$\pi(D)+ \pi(E)=M$, i.e. $\pi(D)=M$ or $\pi(E)=M$ and hence
$D=P$ or $E=P$ as $\Ker{\pi} \ll P$.

Since $M$ is injective, $P$ is divisible, because for all $0\neq x
\in D$ $$\pi(xP)=x\pi(P)=xM=M,$$ i.e. $xP=P$ as $\pi$ has a small
kernel. As $D$ is a principal ideal domain, $P$ is an indecomposable
injective $D$-module and again by Matlis theorem $P\simeq Q$ or
$P\simeq E(D/\qr)$ for some maximal ideal $\qr$. In the later case
we must have $\pr=\qr$ since $$D/\pr=\mathrm{Soc}(E(D/\pr)) \simeq
\mathrm{soc}(E(D/\qr)/\Ker{\pi})=(D/\qr^{i+1}) / (D/\qr^i) \simeq
D/\qr$$ as $E(D/\qr)$ is uniserial and all its submodules are of the
form $D/\qr^i$.
\end{proof}

\subsection{}
The {\it divided power coalgebra} is the path coalgebra $C$ associated to the graph
$$\xygraph{1 :@(ul,ur)^x "1"}.$$
that is the  coalgebra over a field $k$ with basis $\{ 1, x, x^2, \ldots, x^i, \ldots \}$ and comultiplication:
$$ \Delta(x^n)=\sum_{i=0}^n x^i \otimes x^{n-i}$$
and counit
$$ \epsilon(x^n)=\delta_{0,n}.$$

\begin{cor}
Let $C$ be the divided power coalgebra over a field $k$. Then $C$ is its own maximal codense cover in the category of $C$-comodules $\cM$.
\end{cor}
\begin{proof}
The dual algebra $C^* \simeq k[[Z]]$ of $C$ is the power series ring in one variable, by the isomorphism:
$$ f \mapsto \sum_{n=0}^\infty f(x^n)Z^n$$
Note that the power series ring in one variable is a discrete valuation ring, e.g. a principal ideal domain.
Since $C$ is an injective cogenerator in $\cM$ with simple coradical $C_0=k1$,
$C$ is a non-faithful indecomposable injective $C^*$-module over the discrete valuation ring $C^*$.
By Theorem \ref{PID_small_covers}  the only small covers of $C$ in $C^*$-Mod are
$C$ and the quotient field $Q$ of $C^*$. Since $C^*$ is not a $C$-comodule, $Q$ is also not a $C$-comodule. Hence the only small cover of $C$ as $C$-comodule is $C$ itself.
\end{proof}

%

\section{Dual non-singularity of modules}

Recall that a left $R$-module $M$ is called {\it singular} if every element of $M$ is annihilated by an essential left ideal of $R$. An  $R$-module $M$ is called {\it non-singular} if it contains no non-zero singular submodule.

\subsection{}

Non-singularity generalises torsion-freeness of modules to the non-commuta-tive setting.
Lambek's torsion theory is the right concept for a module theoretic setting in which the construction of
maximal dense extension of modules are put. Dual Goldie torsion theories have been studied by various authors
\cite{Ramamurthi}, \cite{generalov}, \cite{Lomp_Spl}.
As singular modules play the r\^ole of torsion modules, small modules will play a similar r\^ole in the dual
situation.
Let $\cS$ be the class of {\it small modules} in $\sigma[M]$, i.e. those which are small in their injective hull in
$\sigma[M]$. $\cS$ is a {\it Serre class}, i.e. it is closed under submodules, factor modules and  extensions (and hence
also under finite direct sums).
Define $$\rho(X)=\Re{X}{\cS} = \bigcap \{ U\subseteq X \mid X/U \in \cS \}$$
for any $X\in \sigma[M]$ and call $X$ {\it dual non-$M$-singular} if $\rho(X)=X$.
These are precisely those modules which do not have any non-zero small homomorphic image.

Since an injective module is a direct summand in any extension, injectives are never small.
Hence {\it cohereditary}\label{def_cohereditary} modules, i.e. those all whose factor modules are injective,
are examples of dual non-$M$-singular modules.
On the other hand there exist examples of injective modules that
are subdirect products of their $M$-small factor modules (see Zoeschinger \cite{zoschinger}).

\subsection{}
Pushing singularity to smaller categories like $\sigma[M]$ needed a
characterisation that was free of refeering to left ideals of a
ring. Concepts for Singularity and their duals had been already
proposed in some abelian categories by Pareigis \cite{Pareigis} and
it is not difficult to see that in the module case a module $M$ is
singular if and only if it is a factor module  of a module by an
essential submodule. In the case of $\sigma[M]$ it turned out, as
shown in \cite{wisbauer_Bimodule}, that non-singularity of $M$ could
be characterised by the internal property that any essential
submodule is dense. This property has been studied by
Zelmanowitz in \cite{Zelmanowitz} where he also termed it {\it
polyform}. It is not difficult to dualise those notions, but it
turns out that they are not always equivalent.

\subsection{}
Dual to a polyform module, call a module $M$ {\it copolyform} if for every small submodule $K$ of $M$, the canonical projection $M\rightarrow M/K$ is a codense cover.
Note that dual non-$M$-singular modules $X$ in $\sigma[M]$
are copolyform since for any small submodule $K$ of $X$ any factor module $K/L$ is also $M$-small
and thus $\Hom{X}{K/L}=0$, i.e. $X\rightarrow X/K$ is codense.
The converse is not true, e.g. $\ZZ$ is copolyform, but not non-$\ZZ$-small.
Copolyform modules had been introduced in \cite{Lomp} and were studied also in \cite{TalebiVanaja}.

\subsection{}
By definition it is clear that copolyform modules can be characterised by their homomorphisms to factor modules.
For any two modules $X$ and $Y$ set
$$\nabla(X,Y)=\{f\in \Hom{X}{Y} \mid \Im{f}\ll Y\}.$$
This set has been introduced by Beidar and Kasch in
\cite{BeidarKasch} were it was termed the {\it cosingular ideal} of $X$
and $Y$. Suppose  $M$ is copolyform and $f\in \nabla(M,M/N)$ for
some $N\ll M$ then $\Im{f}=K/N \ll M/N$ and $N\ll M$ implies $K\ll
M$. But as the projection $M\epim M/K$ is codense,
$f\in\Hom{M}{K/N}=0$. Thus $\nabla(M,M/N)=0$. On the contrary, if
$\nabla(M,M/N)=0$ for all $N\ll M$ then  for any  small cover $\pi:M\epim
F$ with $K=\Ke \pi \ll M$ and submodule $L\subseteq K$ we have $\Hom{M}{K/L}\subseteq \nabla(M,M/L)=0$. Hence
$\pi:M\epim F$ is a codense cover. We have just proved the following
statement:

\begin{thm}
An $R$-module $M$ is copolyform if and only if $\nabla(M,M/N)=0$ for all $N\ll M$.
\end{thm}

Choosing $N=0$ in the above Theorem, we get that a copolyform module
has no non-zero homomorphism with small image, i.e. $\nabla(M):=\nabla(M,M)=0$.
Note that under some suitable projectivity conditions $\nabla(M)$
equals  $\Jac{\End{M}}$.

\subsection{}
Note that for self-projective modules $M$, $\nabla(M)=\Jac{\End{M}}$ (see \cite{wisbauer}). 
\begin{thm}\label{selfprojective_copolyform}
A self-projective module $M$ is copolyform if and only if $\Jac{\End{M}} = 0$.
\end{thm}

Thus a ring $R$ is copolyform as left $R$-module if and only if it is semiprimitive.

\subsection{}\label{facts_on_coalgebras}
Since our aim is to apply the module theoretic terms above to the situation of coalgebras, recall that any coalgebra $C$ of a field $k$ is an injective cogenerator in the category $\cM$ of right $C$-comodules.
Moreover there exists an anti-isomorphism of rings between the dual algebra $C^*$ and the endomorphism of
$C$ as right $C$-comodule and an isomorphism of rings between $C^*$ and the endomorphism of $C$ as left $C$-comodule:

$$\End{_{C^*}C}^{op}\simeq C^* \simeq \End{C_{C^*}}.$$

Under some light injectivity  and cogenerator properties we can say much more about copolyform modules. A module $Q$ is called {\em pseudo-injective} with respect to
a non-zero monomorphism $f:Y\mono X$ if for all non-zero
$g:Y\rightarrow Q$ there exist $h\in \End{Q}$ and $k\in
\Hom{X}{Q}$ such that $fk=gh\neq 0$. A module $Q$ is called {\em
pseudo-injective} in $\sigma[M]$ if it is pseudo-injective with
respect to all non-zero monomorphism $f:Y\mono X$ in $\sigma[M]$.

\begin{lem}\label{Lemma_pseudo-injective_nabla} Let $M$ be pseudo-injective in $\sigma[M]$. Then $\Hom{M/N}{M} = 0$ for all submodules $N$ such that $M/N$ is $M$-small provided $\nabla(M)=0$.
\end{lem}

\begin{proof}
Assume that $M/N$ is small in some module $X\in \sigma[M]$ and let
$f:M/N \rightarrow M$ be a homomorphism. Suppose $f$ is non-zero
then by pseudo-injectivity there are homomorphisms $h\in \End{M}$
and $k\in\Hom{X}{M}$ such that $fh=ik\neq 0$ where $i:M/N\mono X$
denotes the inclusion. Since homomorphic images of small modules are
small, $\Im{fh}=\Im{ik}\ll M$. Considering the projection $p:M\epim
M/N$ we get a homomorphism $pfh\in \End{M}$ whose image is small in
$M$. Since $\nabla(M)=0$, $pfh=0$ which implies $fh=0$, a
contradiction. Thus $\Hom{M/N}{M}=0$.
\end{proof}

\subsection{}
Lemma \ref{Lemma_pseudo-injective_nabla} shows that a
pseudo-injective module $M$ with  $\Hom{M/N}{M} \neq 0$ for all
non-zero $N\subseteq M$, is dual non-$M$-singular if and only if
$\nabla(M)=0$. We will show that this is also equivalent to
$\End{M}$ being non-singular. Say that a module $M$ is  {\em
coretractable} if for all non-zero submodules $N$ of $M$:
$\Hom{M/N}{M}\neq 0$. We first need the following Lemma

\begin{lem}\label{Lemma_non_singular} Let $M$ and $Q$ be left $R$-modules and $T:=\End{Q}$.
Denote by $Z(M^*)$ the singular submodule of $M^*:=\Hom{M}{Q}$ as right $T$-module. Suppose that $Q$ is coretractable then $$Z(M^*)\subseteq \nabla(M,Q)$$
holds. If moreover $Q$ is pseudo-injective with respect to all monomorphisms of the form $g: Q/\Ke g \mono Q$ for any $0\neq g\in T$ then equality hold, i.e. $Z(M^*) = \nabla(M,Q).$
\end{lem}

\begin{proof}
Take $f\in Z(M^*)$. Then $Ann_T(f)=\{g\in T \mid fg=0\}$ is
essential in $T$. Suppose $\Im{f}+U=Q$ for some submodule $U$ of
$Q$. Then $Ann_T(f)\cap Ann_T(U)=Ann_T(\Im{f}+U)=0$ implies
$\Hom{Q/U}{Q}=Ann_T(U)=0$. By hypothesis $U=Q$, i.e. $\Im{f}\ll Q$
and $f\in \nabla(M,Q)$.

\smallskip

Now assume that $Q$ is pseudo-injective with respect to all monomorphisms $g:Q/\Ker{g} \mono Q$.
Let $f\in \nabla(M,Q)$ and $g\in T$ such that $gT\cap
Ann_T(f)=0$. Suppose there exists a non-zero $h\in Ann_T(\Ke
g)\cap Ann_T(\Im f)$. As $h$ defines a non-zero homomorphism from
$Q/\Ke g$ to $Q$ we have by hypothesis endomorphisms $k,l \in T$
such that $0\neq gk=hl$. But as $h\in Ann_T(\Ke g)\cap Ann_T(f)$,
we have $hl=gk\in gT\cap Ann_T(f)=0$; a contradiction. Thus
$Ann_T(\Ke g)\cap Ann_T(\Im f)=0$ and
$$0 = Ann_T(\Ke g)\cap Ann_T(\Im f) =  Ann_T(\Ke g + \Im f) \simeq \Hom{Q/(\Ke g + \Im f)}{Q}.$$
 Since $Q$ is coretractable, $\Ke{g}+\Im{f}=Q$, but as $\Im f \ll Q$, $g=0$.
\end{proof}

Note that the condition in Lemma \ref{Lemma_non_singular}(2) is
fulfilled if $Q$ is semi-injective, i.e. injective with respect to
all monomorphisms of the above form, or if $Q$ is pseudo-injective
in $\sigma[Q]$.

\subsection{}
The last Lemma \ref{Lemma_non_singular} together  with \ref{Lemma_pseudo-injective_nabla} enables us to characterise those copolyform modules which are injective cogenerators:

\begin{thm}\label{copolyform_injcog} Let $M$ be a coretractable left $R$-module that is
pseudo-injective in $\sigma[M]$. Then the following statements are
equivalent:
\begin{enumerate}
\item[(a)] $M$ is dual non-singular in $\sigma[M]$.
\item[(b)] $M$ is copolyform.
\item[(c)] $\nabla(M)=0$.
\item[(d)] $\End{M}$ is a right non-singular ring.
\end{enumerate}
\end{thm}

\subsection{}
The lattice of submodules of a module is pseudo-complemented, but
its dual lattice does not need to be. To overcome this problem while
dualising module theoretic notions, one has to make suitable
assumption on the lattice of submodules. An $R$-module $M$ is called
{\it weakly supplemented} if any submodule $N$ of $M$ has a {\it
weak supplement}, that is a submodule $L$ of $M$ such that $N+L=M$
and $N\cap L\ll M$. This is a weak form of a pseudo-complement in
the dual lattice of submodules of $M$.

\begin{thm}
The following statements are equivalent for a weakly supplemented module:
\begin{enumerate}
    \item[(a)] $M$ is copolyform.
    \item[(b)] $\nabla(M,M/N)=0$ for all $N\subseteq M$.
    \item[(c)] Every factor module of $M$ is copolyform.
    \item[(d)] $\nabla(M)=0$ and $M$ is {\it $M$-im-small-projective}, i.e. any diagram
$$\xymatrix{ & M\ar[d]^f\ar@{-->}[dl]^h \\ M \ar[r]_g & L \ar[r] & 0 }$$
with $\Im{f}\ll L$ can be commutatively extended by some $h:M\rightarrow M.$
\end{enumerate}
\end{thm}

\begin{proof}
$(a)\Rightarrow (b)$ Let $f:M\rightarrow M/N$ have small image and choose a weak supplement $L$ of $N$. Thus
$M/N=(N+L)/N \simeq L/(L\cap N) \subseteq M/(L\cap N)$. Since $L\cap N\ll M$ and $f\in\nabla(M,M/(L\cap N)$, we have $f=0$ by $(a)$.

$(b)\Rightarrow (c)$ let $N\subseteq L \subseteq M$ such that $L/N\ll M/N$ and $f\in\nabla(M/N,M/L)$. Then $f\pi_N \in \nabla(M,M/L)=0$, i.e. $f=0$. Hence $M/N$ is copolyform.

$(c)\Rightarrow (a)$ is trivial and
$(b)\Rightarrow (d)$ is clear, since for $N=0$, $\nabla(M,M)=\nabla(M)=0$ and as $\nabla(M,M/N)=0$ for all factor modules $L$ of $M$, there are no non-zero homomorphisms $f:M\rightarrow L$ with small image, i.e. $M$ is trivially $M$-im-small projective.

$(d)\Rightarrow (a)$ Let $f\in \nabla(M,M/N)$ with $N\ll M$ and denote by $\pi_N:M\rightarrow M/N$ the canonical projection. By $M$-im-projectivity there exists $h:M\rightarrow M$ such that $\pi_Nh = f$. Since
$\Im{f}=\Im{\pi_Nh} \ll M/N$ and $N\ll M$, we have $\Im{h}\ll M$, i.e. $h\in \nabla(M)=0$. Thus $f=0$.
\end{proof}

A module which satisfies condition $(c)$ is also called {\it strongly copolyform}. This is in general a stronger condition then copolyformness. In \cite{TalebiVanaja} strongly copolyform modules are called copolyform.


\subsection{}
A module $M$ is called {\it couniform} or {\it hollow}
if $N+L=M$ implies $N=M$ or $L=M$ for all proper submodules $N, L$ of $M$. Uniserial modules are couniform  and couniform modules are  indecomposable. Furthermore couniform modules are trivially weakly supplemented since all proper submodules are small. From the last characterisation of copolyform modules
we easily deduce that a couniform module is copolyform if and only if
every projection $M\epim M/N$ for any proper submodule $N$ of $M$ is codense. Couniform
copolyform modules are called {\it epiform} and satisfy the property that all of
 their non-zero endomorphisms are epimorphisms.
The converse holds under some suitable assumptions as we will see later.

\subsection{}

In case of couniform modules we deduce from \ref{copolyform_injcog}
the following

\begin{cor}\label{epiform_injcog} Let $M$ be a couniform coretractable left $R$-module that is
pseudo-injective in $\sigma[M]$. Then the following statements are
equivalent:
\begin{enumerate}
\item[(a)] $M$ is dual non-$M$-singular.
\item[(b)] $M$ is epiform.
\item[(c)] Every non-zero endomorphism of $M$ is an epimorphism.
\item[(d)] Every non-zero homomorphism from a factor module $L$ of $M$ to $M$ is surjective.
\item[(e)] $\End{M}$ is a domain.
\end{enumerate}
\end{cor}

\begin{proof}
$(a)\Leftrightarrow (b)$ follows from \ref{copolyform_injcog}.\\
$(b)\Rightarrow (c)$ For any $0\neq f\in\End{M}$, $\Im{f}\not\ll M$ as $\nabla(M)=0$. Thus $\Im{f}=M$.\\
$(c)\Rightarrow (d)$ Let $f:M/N\rightarrow M$ since $\pi_Nf$ is an epimorphism of $M$, $f$ has to be an epimorphism (here $\pi_N$ denotes the projection).\\
$(d)\Rightarrow (e)$ If $fg=0$, then $\Im{f}\subseteq \Ker{g}$. And if $f\neq 0$, then $M=\Im{f}=\Ker{g}$, i.e. $g=0$.\\
$(e)\Rightarrow (a)$ follows from \ref{copolyform_injcog} as domains are non-singular.
\end{proof}

\subsection{}
Copolyform module with projective covers can be characterise by their endomorphism rings.

\begin{prop} \label{prop_2} Let $M$ be an $R$-module with projective cover $P$ in $\sigma[M]$.
  Then $M$ is copolyform if and only if  $\Jac{\End{P}}=0$.
\end{prop}

\begin{proof}
Recall that $\Jac{\End{P}} = \nabla(P)$.
Assume $M$ to be copolyform and let $f \in \nabla(P)$.  Then, for
  any $g \in \Hom{P}{M}$, $U:=\Im{fg} \ll M$. However, by Lemma \ref{prelema}, $\Hom{P}{U}=0$ and so $fg=0.$
 This implies $\Im{f} \subseteq \Ker{g}$ and so
  $$\Im{f} \subseteq \bigcap\left\{\Ker{g}:g \in \Hom{P}{M}\right\}  = \Re{P}{M}=0,$$
as $P$ is cogenerated by $M$ (see \cite[18.4]{wisbauer}). Thus $f=0$, i.e.  $\Jac{\End{P}} = 0$.\\
On the contrary if $\nabla(P)=0$, then $P$ is copolyform by \ref{selfprojective_copolyform}.
Denote by $p:P\epim M$ the projection and let $\pi:M\epim X$ be any small cover. The composition $p\pi:P\epim X$ is also
a small cover and therefore codense. In particular $\Hom{P}{\Ker{p\pi}}=0$ and, by projectivity of $P$, $\Hom{P}{\Ker{\pi}}=0$. By \ref{projectiveCoverCodense} $\pi$ is a codense cover, i.e. $M$ is copolyform.
\end{proof}

\subsection{}
The last proposition showed that a projective cover of a copolyform module is copolyform as well.

\begin{cor}
Let $M$ be a copolyform module with projective cover $P$ in $\sigma[M]$, then $\End{M}$ is a subring of $\End{P}$ such that every epimorphism $f\in \End{M}$ with small kernel is invertible in $\End{P}$.
\end{cor}
\begin{proof}
Denote by $p:P\rightarrow M$ the projection and take any non-zero $f\in \End{M}$. Then by the projectivity of $P$, there exists a non-zero ${\bar f} \in \End{P}$ such that $pf=\bar{f}p$. Suppose there exists another $g\in \End{P}$ such that $pf=gp$, then
$$ 0 = pf-pf = (\bar{f}-g)p$$
implies $\Im{\bar{f}-g}\subseteq \Ker{p}$, i.e. $\bar{f}-g \in \nabla(P)=0$. Hence $\bar{f}=g$. Thus the correspondence $f\mapsto \bar{f}$ is uniquely defined.

Now assume that $f$ is an epimorphism with small kernel, then $pf=\bar{f}p$ implies that $\bar{f}p$, and hence $\bar{f}$ is an epimorphism with small kernel. By the projectivity of $P$, $\bar{f}$ splits and, as $\Ker{\bar{f}}\ll P$, must be  an isomorphism.

\end{proof}

\subsection{}

The existence of a projective cover, turns the class $\cS$ of $M$-small modules
into a cotorsion class:

\begin{prop} Assume that $M$ is dual non-singular in $\sigma[M]$ and has a projective cover $P$ in
  $\sigma [M]$. Then the class of small modules in $\sigma[M]$
  is closed under submodules, factor modules, extensions and  direct
products (in $\sigma [M]$) and can be described as:
$$ \cS = \{ X \in \sigma [M] : \Hom{P}{X}=0\}$$
 Moreover for any $Z\in \sigma[M]$,
$\rho(Z) = \Re{Z}{\cS}$ is dual non-$M$-singular and $Z/\rho(Z)$ is
$M$-small.
  \end{prop}

\begin{proof}
Note that if $\cS$ can be described as stated above, then it also
satisfies the closure properties. Hence we only need to show that
$\cS$ equals the class of modules $X$ with $\Hom{P}{X}=0$.  Let $X$ be any
module in $\sigma[M]$ and $\widehat{X}$ its injective hull in
$\sigma[M]$. By \cite[17.9]{wisbauer}, ${\widehat X}$ is
$M$-generated and hence $P$-generated. If $X$ is not $M$-small, then
it is not small in its $M$-injective hull $\widehat X$. Thus assume
there is a proper submodule $Y$ of ${\widehat X}$ such that
$X+Y={\widehat X}$. Then $X/(X\cap Y) \simeq {\widehat X}/Y$ is a
nonzero $P$-generated $R$-module. Hence there is an index set
$\Lambda$ and an epimorphism $f:P^{(\Lambda )}\rightarrow
X/(X \cap Y)$ and so, since $P^{(\Lambda )}$ is projective in
$\sigma[M]$, $f$ can be lifted to a homomorphism $g:P^{(\Lambda)}
\rightarrow X$, i.e. $\Hom{P}{X}\neq 0$. Hence $X\not\in\cS
\Rightarrow \Hom{P}{X}\neq 0$.

On the other hand assume $0\neq X\in \cS$ and $f\in \Hom{P}{X}$. Denote
by $Y=\Im{f}$ and let $\pi:P\rightarrow M$ be the projection. Then
extend $f$ to a homomorphism
$$ g: M \simeq P/\Ker{\pi}\rightarrow Y/(\Ker{\pi})f$$
sending $p+\Ker{\pi}$ to $(p)f+(\Ker{\pi})f$. Since $M$ is dual
non-M-singular, $g=0$ and $Y=\Im{f}\subseteq (\Ker{\pi})f$. Thus
$P=\Ker{\pi}+\Ker{f}$, but since $\Ker{\pi}\ll P$, $\Ker{f}=P$ and
$f=0$. This shows that $X\in\cS$ implies $\Hom{P}{X}=0$ proving the
equality of the  classes indicated.

Thus $\cS$ is closed under submodules, factor modules, direct
products and extensions. Note that it follows also that $P$ is dual
non-$M$-singular. Moreover since $Z/\rho(Z)$ is a subdirect product
of $M$-small modules, it is $M$-small. Furthermore, since $P$ is
projective and $$\Hom{P}{\rho(Z)/\Tr{P}{\rho(Z)}}=0,$$  we must have
$\rho(Z)=\Tr{P}{\rho(Z)}$, i.e. $\rho(Z)$ is $P$-generated and
therefore dual non-$M$-singular.
\end{proof}

In the case above, $P$ generates the {\it cotorsion theory} whose
cotorsion modules are the $M$-small modules in $\sigma[M]$. the
cotorsion free modules are precisely the $P$-generated modules.

\section{Non-singular Coalgebras}

Having defined a dual non-singularity concept for modules, we are going to apply it to comodules.
Let $C$ be a coalgebra over a field $k$. Any right $C$-comodule $M$ carries a natural left $C^*$-module structure.
Call a right(left) $C$-comodule $M$ copolyform (resp. epiform) if it is copolyform (resp. epiform) as left(right)
 $C^*$-module.

Theorem \ref{copolyform_injcog} and the facts on coalgebras \ref{facts_on_coalgebras} yield  the following

\begin{thm} Let $C$ be a coalgebra over a field $k$. Then the following statements are equivalent:
\begin{enumerate}
\item[(a)] $C$ is a copolyform right $C$-comodule.
\item[(b)] $\End{_{C^*}C}$ is a right non-singular ring.
\item[(c)] $C^*$ is a left non-singular ring.
\item[(d)] $\End{C_{C^*}}$ is a left non-singular ring.
\item[(e)] $C$ is a copolyform left $C$-comodule.
\end{enumerate}
Any coalgebra that satisfies one of the above conditions is called {\it non-singular}.
\end{thm}

\begin{proof}
$(a)\Leftrightarrow (b)$ follows from Theorem \ref{copolyform_injcog}.\\
$(b)\Leftrightarrow (c)$ follows from the anti-isomorphism between $\End{_{C^*}C}$ and $C^*$.\\
$(c)\Leftrightarrow (d)$ follows from the isomorphism between $\End{C_{C^*}}$ and $C^*$.\\
$(d)\Leftrightarrow (e)$ follows from Theorem \ref{copolyform_injcog} (for right $R$-modules).
\end{proof}

\subsection{}
In \cite{NastasescuTorrecillasZhang}, Nastasescu, Torrecillas and Zhang called a coalgebra $C$
{\it hereditary} if $C$ is a cohereditary left (and/or right)  $C$-comodule.
By our remark in \ref{def_cohereditary} cohereditary modules are
dual non-singular. Hence by \ref{copolyform_injcog} any hereditary coalgebra is non-singular.
Chin showed in \cite{Chin} that any path coalgbera is hereditary.
Furthermore Chin and Montgomery showed in \cite{ChinMontgomery} that any coalgebra over an algebraically closed field is Morita-Takeuchi
equivalent to a subcoalgebra of a path coalgebra. Thus hereditary and hence non-singular coalgebras are ubiquitous.

\subsection{}
In \cite{NastasescuTorrecillasZhang} it has been also proven that a finite dimensional coalgebra $C$ is
hereditary if and only if $C^*$ is left hereditary.
Since there are finite dimensional algebras which are left non-singular, but not left hereditary, we can construct
coalgebras which are non-singular but not hereditary. Let $k$ be a field let $R$ be any finite dimensional $k$-algebra
which is not left hereditary; for example $R=k[x]/(x^2)$. Then define
$$ A=\left( \begin{array}{cc} k & 0 \\ R & R  \end{array} \right).$$
By \cite[4.4.3]{goodearl}, $A$ is right non-singular, but not right hereditary  by \cite[4.4.7]{goodearl} as $R$ is not right hereditary.
Hence $C = \left( A^{op} \right)^*$ is a non-singular coalgebra which is not hereditary.

\subsection{}

Call a coalgebra $C$ {\it cosemiprime} if $I\wedge I \neq C$ holds for all proper subcoalgebras $I$ of $C$.
It is not difficult to see that $C$ is a cosemiprime coalgebra if and only if $C^*$ is semiprime and we deduce that a cocommutative coalgebra is non-singular if and only if $C$ is cosemiprime.

\subsection{} The strict hierarchie of coalgebraic properties
\begin{center}
\textit{cosemisimple} $\Rightarrow $ \textit{hereditary} $\Rightarrow $ \textit{non-singular}
\end{center}
collapses when assuming some flatness  condition on the coalgebra:
Since a coalgebra $C$ is flat as right $C^*$-module if and only if $C^*$ is left self-injective
(see \cite{BrzezinskiWisbauer}), we have that the dual algebra $C^*$ of
a non-singular coalgebra $C$ which is flat as left $C^*$-module must be a
left self-injective and  left non-singular ring and hence von Neumann regular
(as it equals its own maximal left ring of quotient).
Note that a von Neumann regular ring is semiprimitive, hence $\Jac{C^*}=0$.
By \cite{BrzezinskiWisbauer}, $\Jac{C^*}=C_0^\perp$ where $C_0$ denotes the coradical of $C$.
Hence $C_0^\perp=0$ implies $C=C_0$.
We just proved the following theorem:
\begin{thm}\label{vonNeumannCoalgebras}
A coalgebra $C$ is cosemisimple if and only if $C$ is non-singular and flat as right $C^*$-module.
\end{thm}

Since finite dimensional Hopf algebras are projective as comodule, we deduce
that finite dimensional Hopf algebra which are right non-singular coalgebras are cosemismple.

\subsection{}

The characterisation \ref{epiform_injcog} of epiform modules yields that a coalgebra $C$ is
epiform as right (or left) comodule if and only if $C^*$ is a domain.
Recall that a coalgebra $C$ is called {\it coprime} if $C^*$ is a prime ring.
As we see, any coalgebra that is epiform as coalgebra is a coprime coalgebra.
In case $C$ is cocommutative those notions are equivalent.

\subsection{}

We are going to show that there exists a dichotomie for coprime coalgebras that states that over a coprime coalgebra either every comodule is projective or no non-trivial comodule is projective.

\begin{thm}\label{finite}
The following statements are equivalent for a coprime coalgebra $C$ over a field $k$:
\begin{enumerate}
\item[(a)] $C$ is a matrix $k$-coalgebra, i.e. $C^*$ is a matrix algebra over a division ring.
\item[(b)] $C^*$ is a simple ring.
\item[(c)] $C$ is finite dimensional.
\item[(d)] Every non-zero right (left) $C$-comodule is projective as $C^*$-module.
\item[(e)] There exists a non-zero projective right (left) $C$-comodule.
\item[(f)] No non-zero right (left) $C$-comodule is singular as $C^*$-module.
\item[(g)] There exists a non-zero right (left) $C$-comodule that is not singular as a $C^*$-module.
\item[(h)] Every right (left) $C$-comodule is injective.
\end{enumerate}
\end{thm}
\begin{proof}
(a) $\Rightarrow$ (b) is clear.

(b) $\Rightarrow$ (c) assume $C^*$ is simple, then $C$ is a simple coalgebra, because if $D$ is any subcoalgebra of  $C$, then $D^\perp$ is an ideal of $C^*$ and hence $0$ or $C^*$ ,i.e. $D=C$ or $D=0$. Since any non-zero
element of $C$ is contained in a non-zero finite dimensional subcoalgebra of $C$, $C$ must be finite dimensional.

(c) $\Rightarrow$ (a) since $C$ is finite dimensional, $C^*$ is finite dimensional. As $C^*$ is also a prime ring, it must be a matrix algebra.

(a) $\Rightarrow$ (d) is clear.

(d) $\Rightarrow$ (e) $\Rightarrow$ (g) and 
(d) $\Rightarrow$ (f) $\Rightarrow $ (g) are trivial since projective modules are not singular.

(g) $\Rightarrow$ (c) Suppose $M$ is a non-zero left $C$-comodule which is not singular as $C^*$-module.
Then there exists a $C^*$-submodule $N$ of $M$ which is not singular. We might choose $N$ to be a cyclic
$C^*$-submodule of $M$. Since comodules are locally finite dimensional, $N$ is finite dimensional.
As the annihilator $\Ann_{C^*}(N)$ is not an essential left ideal of $C^*$, but all non-zero ideals of a
prime ring are essential as left ideals, we conclude that $\Ann_{C^*}(N)=0$, thus
$$C^* = C^*/\Ann_{C^*}(N) \hookrightarrow \oplus_{i=1}^{s} C^*/ Ann_{C^*}(n_i)$$
is finite dimensional, where ${n_i}$ is a generating set of $N$.

$(a)\Leftrightarrow (h)$ is clear, since $C$ is cosemisimple.

\end{proof}


\subsection{}\label{Infinite_coprime_coalgebras}
By negating $(c),(e)$ and $(g)$ we get of the last Theorem we deduce the following

\begin{cor}
The following statements are equivalent for a coprime coalgebra $C$ over a field $k$.
\begin{enumerate}
\item[(a)] $C$ has infinite dimension.
\item[(b)] Every right or left $C$-comodule is singular as $C^*$-module.
\item[(c)] There is no non-zero projective object in the category of right or left $C$-comodules.
\end{enumerate}
\end{cor}

\subsection{}
The last Corollary shows the dichotomie of coprime coalgebras: Either every comodule is coalgebra and the coalgebra is necessarily a matrix coalgebra or every comodule is singular as $C^*$-module and $\cM$ has no non-zero projective object.

This dichotomie shows also that we can not use projective cover to build maximal codense covers of infinite dimensional coprime coalgebras.

\subsection{}
From \ref{epiform_injcog} have that any  $C$ which is epiform is either the dual of a
finite dimensional divison algebra $K$ over $k$ or infinite dimensional such that the category of right $C$-comodules consists
of torsion $C^*$-modules, in particular there are no non-zero projective objects in $\cM$.

\subsection{}

Note that any coalgebra $C$ can be written as a sum of indecomposable injective comodules $E_\lambda$.
If $C$ is cocommutative then each of the $E_\lambda$ is actually a subcoalgebra of $C$.
Assume now that $C$ is a cocommutative semiperfect coalgebra over a field $k$, then  $C=\bigoplus_\lambda E_\lambda$ is a direct coproduct of finite dimensional cocommutative indecomposable coalgebras. If moreover $C$ is non-singular, then each of the $E_\lambda$ is also non-singular and $E_\lambda^*$ is a finite dimensional commutative semiprime $k$-algebra. Thus $E_\lambda^*$ is a finite field extension $K_\lambda$ of $k$ and $E_\lambda = K_\lambda^*$ is a finite dimensional simple coalgebra. Thus we have proved the following

\begin{thm}
Any cocommutative non-singular and semiperfect coalgebra is cosemisimple.
\end{thm}

\end{document}